\documentclass[11pt]{article}
\usepackage{amsthm,amsfonts,amsmath,amssymb}
\usepackage{mdwlist}
\usepackage{mathrsfs,tikz}
\usepackage{dsfont}
\newtheorem{thm}{Theorem}[section]

\newtheorem{lem}[thm]{Lemma}
\newtheorem{cor}[thm]{Corollary}
\newtheorem{defi}[thm]{Definition}

\newtheorem{conj}[thm]{Conjecture}

\numberwithin{equation}{section}

\def\demo{\noindent{\bf Proof}\hskip10pt}

\def\qed{\hfill $\Box$}
\def\lg{\langle}
\def\rg{\rangle}
\def\rr#1{\item[{\rm (#1)}]}

\def\Ker{\hbox{\rm Ker}}
\def\Im{\hbox{\rm Im}}

\begin{document}
\title{On quadratic conjecture\thanks{This work was supported by NSFC (No. 11971280 \& 11771258)}}
\author{Jingjing Duan and Lijian An\thanks{Corresponding author. e-mail: anlj@sxnu.edu.cn }\\
Department of Mathematics,
Shanxi Normal University\\
Taiyuan, Shanxi 030031, P. R. China\\
 }

\maketitle

\begin{abstract}

Quadratic conjecture is a strengthening of oliver's $p$-group conjecture. Let $G$ be a $p$-group of maximal class of order $p^n$. We prove that if $n\le 8$ or $n\ge \max\{2p-6,p+2\}$ then $G$ satisfies
Quadratic Conjecture. Hence quadratic conjecture holds if $G$ is a $p$-group of maximal class where $p\le 7$.


\end{abstract}

\baselineskip=16pt

\section{Introduction}

 An open question in the theory of $p$-local finite groups claims that every fusion system has a unique $p$-completed classifying space. In \cite{O}, Oliver introduced the characteristic subgroup $\mathfrak{X}(S)$ for a finite $p$-group $S$. For odd primes he demonstrated that the conjecture outlined below would substantiate the unique existence of the classifying space. Recall that $J(S)$ denotes the Thompson subgroup of $S$ generated by all elementary abelian subgroups of the greatest rank.

\medskip

\begin{conj}{\rm \cite[Conjecture 3.9]{O}}\label{oliver} For any odd prime $p$ and any $p$-group $S$, $J(S)\le \mathfrak{X}(S)$.
\end{conj}

Conjecture \ref{oliver} is called Oliver's p-group conjecture.
Let $G$ be a finite $p$-group and $V$ an elementary abelian $p$-group
on which $G$ acts faithfully. Then $V$ is an {\it F-module} for $G$ if there exists a non-trivial
elementary abelian subgroup $E$ of $G$ such that $|E|\cdot |C_V (E)| \ge |V |$. In this case, we call $E$
an {\it offender}.

In \cite{GHL}, Green, H$\acute{e}$thelyi and Lilienthal obtained the following reformulation of Conjecture \ref{oliver}.

\begin{conj}{\rm\cite[Conjecture 1.3]{GHL}}\label{Oliver conjecture} Let $p$ be an odd prime and $G$ a finite $p$-group. If the faithful $\mathbb{F}_p[G]$-module $V$ is an $F$-module, then
there is an element $1\ne g\in\Omega_1(Z(G))$ such that the minimal polynomial of the action of $g$ on $V$ divides $(x-1)^{p-1}$.
\end{conj}

An element $g\in G$ is said to be
{\it quadratic} on the $\mathbb{F}_p[G]$-module $V$ if its action has minimal polynomial $(X-1)^2$. Note that if $V$ is
faithful then quadratic elements must have order $p$.

In \cite{GHM2}, Green, H$\acute{e}$thelyi and Mazza propose the following strengthening of Conjecture \ref{Oliver conjecture}.
%\begin{rem} F-module is short for ¡°failure of (Thompson) factorization module¡±.
%\end{rem}
\medskip

\noindent {\bf Quadratic conjecture}\ {\rm\cite[Conjecture 1.4]{GHM2}}\label{quadratic conjecture}\ Let $p$ be a prime and $G$ a finite $p$-group. If the faithful $\mathbb{F}_p[G]$-module $V$ is an $F$-module, then
there are quadratic elements in $\Omega_1(Z(G))$.

\medskip

Observe that Quardatic Conjecture is trivially true for $p=2$. For $p=3$, Quardatic Conjecture is just Conjecture \ref{Oliver conjecture}.  By \cite[Theorem 1.2]{GHM}, Quadratic Conjecture holds if $G$ is metabelian or of (nilpotence) class at most four. By \cite[Theorem 1.3]{GHM}, Conjecture \ref{Oliver conjecture} holds if $G$ is a $p$-group of maximal class. Hence Quadratic Conjecture holds if $G$ is a $3$-group of maximal class. In this paper, we get the following main result:

\begin{thm}\label{main} Let $p$ be an odd prime and $G$ a $p$-group of maximal class of order $p^n$. If $n\le 8$ or $n\ge \max\{2p-6,p+2\}$ then $G$ satisfies
Quadratic Conjecture.
\end{thm}

By Theorem \ref{main}, quadratic conjecture holds if $G$ is a $p$-group of maximal class where $p\le 7$.

\section{Quadratic Offenders}

\begin{defi}{\rm \cite[11, 2.3]{MS}} Let $G$ be a finite group and $V$ a faithful $\mathbb{F}_p[G]$-module. For a
subgroup $H\le G$ one sets
$$j_H(V):=\frac{|H||C_V(H)|}{|V|}\in \mathbb{Q}.$$
\end{defi}
Note that $j_1(V) = 1$. Suppose that $E$ is an element abelian subgroup of $G$. Then $E$ is an offender on $V$ if and only if $j_E(V)\ge 1$.

\begin{defi}{\rm\cite[26.5]{GLS}} Let $G$ be a finite $p$-group and $V$ a faithful $\mathbb{F}_p[G]$-module. We denote by $\mathscr{E}(G)$
the poset of non-trivial elementary abelian subgroups of $G$. Define
$$\mathscr{P}(G, V):= \{E \le G \mid E\in\mathscr{E}(G) \ \mbox{and}\ j_E(V)\ge j_F(V)\  \forall \ 1 \le F\le E\}.$$
\end{defi}

Note that $V$ is an F-module if and only if $\mathscr{P}(G, V)$ is non-empty. The subgroups
in $\mathscr{P}(G, V)$ are sometimes called {\it best offenders}.

A {\it quadratic offender} is a subgroup which is both quadratic
($[V, E, E] = 0$) and an offender.

\begin{thm}{\rm (Timmesfeld's replacement theorem)\cite[Theorem 4.1]{GHM}} Let $V$ be a faithful $\mathbb{F}_p[G]$-module
and $E\in\mathscr{P}(G, V)$. Then there is a quadratic offender $F\in\mathscr{P}(G, V)$ which satisfies $j_F(V) = j_E(V)$ and $F\le E$ .
\end{thm}

By Timmesfeld's replacement theorem, A minimal offender must be a quadratic offender.

Following from \cite{GHM2}, an abelian subgroup $A$ of $G$ is said to be {\it weakly closed} if
$[A, A^g]\ne 1$ holds for every $G$-conjugate $A^g\ne A$. The following theorem show that there is also a weakly closed quadratic offender on F-module $V$.

\begin{thm}{\rm \cite[Proposition 4.5]{GHM2}}\label{weakly}
Suppose that the faithful $\mathbb{F}_p[G]$-module $V$ is an F-module. Set
$$j_0 =\max\{ j_E(V)\mid E \ \mbox{an offender}\}.$$
Then there is a weakly closed quadratic offender $E$ with $j_E(V) = j_0$.
Moreover if $D\le G$ is any offender with $j_D(V)=j_0$, then there is such an $E$ which is a subgroup of the
normal closure of $D$.
\end{thm}

\begin{lem}\label{p3} Let $G$ a finite $p$-group and $V$ an $F$-module. Suppose that there is no quadratic element in $\Omega_1(Z(G))$. Then
 \begin{itemize}
   \rr1 $|E|\ge p^2$ for any offender $E$ on $V$;
   \rr2 $|E|\ge p^3$ for any weakly closed offender $E$ on $V$.
 \end{itemize}
\end{lem}
\demo
(1) Otherwise, there is an offender $E$ of order $p$. Since $V$ is faithful, $C_V(E)\ne V$. Since $E$ is an offender, $|E||C_V(E)|\ge |V|$. Hence $C_V(E)$ is maximal in $V$. Since $C_V(E)$ is $Z(G)$-invariant, $[V,Z(G)]\le C_V(E)$. Since $[V,Z(G)]\le C_V(E^g)$ for all $g\in G$, $[V,Z(G)]\le C_V(E^G)$. Take $1\ne z\in Z(G)\cap E^G$. Then $[V,z,z]=0$. Hence $z$ is a quadratic element in $\Omega_1(Z(G))$. This contradicts the hypothesis.

(2) Otherwise, there is a weakly closed offender $E$ of order $p^2$. By (1), $E$ is a minimal offender. Hence $E$ is also a quadratic offender. If $E\trianglelefteq G$, then $E\cap Z(G)\ne 1$. Hence there is a quadratic element in $\Omega_1(Z(G))$. This contradicts the hypothesis. In the following, we may assume that $N_G(E)\ne G$. Let $g\in N_G(N_G(E))\setminus N_G(E)$. Then $N_G(E^g)=N_G(E)^g=N_G(E)$ and $E^g\ne E$. Hence $E$ and $E^g$ normalize each other. It follows that $[E,E^g]\le E\cap E^g$. Since $E$ is weakly closed, $[E,E^g]\ne 1$ and hence $E\cap E^g\ne 1$.

Let $F=E\cap E^g$. By (1), $F$ is not an offender. Hence $|C_V(E)|\le |C_V(F)|\le \frac{1}{p^2}|V|$. On the other hand, since $E$ is an offender, $|C_V(E)|\ge \frac{1}{p^2}|V|$. Hence $|C_V(E)|=\frac{1}{p^2}|V|$ and $C_V(F)=C_V(E)$. The same reason gives that $C_V(F)=C_V(E^g)$.
Since $E$ and $E^g$ are both quadratic, $[V,E]\le C_V(E)=C_V(F)=C_V(E^g)$ and $[V,E^g]\le C_V(E^g)=C_V(F)=C_V(E)$. It follows that $[V,E,E^g]=[V,E^g,E]=0$. By Three Subgroups Lemma, $[V,[E,E^g]]=0$. Hence $V$ is not faithful. This is a contradiction. \qed

\section{Quadratic elements}

For finding quadratic element, Lemma 4.1 of \cite{GHL} is a key result. In the following, we refine this lemma, and give an elementary proof.

\begin{lem}\label{quadratic element}
Suppose that $p$ is an odd prime, that $G$ is a non-trivial $p$-group, and that
$V$ is a faithful $\mathbb{F}_p[G]$-module. Suppose that $a, b\in G$ are such that $c:=[a,b]$ is a non-trivial element of
$C_G (a,b)$. If $a$ is quadratic, then $[V,E,E]=0$ where $E=\lg a,c\rg$.
\end{lem}
\demo Let $d=a^b=ac$ and $e=a^{b^2}=ac^2=a^{-1}d^2$. Then $E=\lg a\rg\times \lg d\rg$. Since $a$ is quadratic, $d$ and $e$ are also quadratic. Since $[a,d]=1$, $[v,a,d]=[v,d,a]$ for all $v\in V$. Let $v\in V$ and $i,j,k,l\in \mathbb{Z}$. Then
$$[v,a^id^j,a^kd^l]=(il+jk)[v,a,d].$$
In particular, $[v,e,e]=-4[v,a,d]$. It follows from $[v,e,e]=0$ that $[v,a,d]=0$. Hence $[v,a^id^j,a^kd^l]=0$ for all $v\in V$ and $i,j,k,l\in \mathbb{Z}$. Thus we have $[V,E,E]=0$.
\qed

\medskip

Recall that the following terminology. Given a finite group $G$, the ascending central series is
defined inductively by $Z_0(G)=1$ and $Z_{r+1}(G)$ is the normal subgroup of $G$ containing $Z_r(G)$ and such that $Z_{r+1}(G)/Z_r(G)=Z(G/Z_r(G))$, for all $r\ge 1$. The class $c$ of $G$ is the smallest number such that $Z_c(G)=G$. In this paper we use
$G>G_2=G'>\cdots>G_{c+1} = 1$
to denote the descending
central series of $G$, where $G_{r+1}=[G_r,G]$ for all $r\ge 2$.

\begin{lem}\label{Z_i(G)} Let $G$ be a $p$-group and $V$ an F-module. $E$ is a weakly closed quadratic offender. If there is no quadratic element in $\Omega_1(Z(G))$, then
\begin{itemize}
  \rr1 there is no quadratic element in $Z_2(G)$;
  \rr2 $Z_2(G)\cap E=1$, $[Z_2(G),E]=1$;
  \rr3 $[Z_3(G),E]=1$;
  \rr4 $[Z_4(G),E]\le Z_3(G)\cap E$.
\end{itemize}
\end{lem}
\demo (1) Otherwise, let $a$ be a quadratic element in $Z_2(G)$. Then $a\not\in Z(G)$ since there is no quadratic element in $Z(G)$. Hence there is $b\in G$ such that $c:=[a,b]\ne 1$. Since $c\in Z(G)\le C_G(a,b)$, by Lemma \ref{quadratic element}, $c$ is also a quadratic element. This is a contradiction.

(2) By (1), $Z_2(G)\cap E=1$. Since $[Z_2(G),E]\le Z(G)$, $[E^g, E]=1$ for all $g\in Z_2(G)$. Since $E$ weakly closed, $E^g=E$ for all $g\in Z_2(G)$. It follows that $[Z_2(G),E]\le Z(G)\cap E=1$.

(3) Since $[Z_3(G),E]\le Z_2(G)$, by (2), $[Z_3(G),E,E]=1$. It follows that $[E^g, E]=1$ for all $g\in Z_3(G)$. Since $E$ weakly closed, $E^g=E$ for all $g\in Z_3(G)$. Furthermore, $[Z_3(G),E]\le Z_2(G)\cap E=1$.

(4) Since $[Z_4(G),E]\le Z_3(G)$, by (3), $[Z_4(G),E,E]=1$. It follows that $[E^g, E]=1$ for all $g\in Z_4(G)$. Since $E$ weakly closed, $E^g=E$ for all $g\in Z_4(G)$. Furthermore, $[Z_4(G),E]\le Z_3(G)\cap E$.
\qed
\medskip

\begin{cor}{\rm \cite[Theorem 5.2]{GHM}}\label{class 4}
Suppose that $G$ is a $p$-group and $V$ is a faithful $\mathbb{F}_p[G]$-module such that
$\Omega_1(Z(G))$ has no quadratic elements. If $G$ has class at most four, then $V$ cannot be an F-module.
\end{cor}
\demo Otherwise, $V$ is an F-module. By Lemma \ref{Z_i(G)} (4), $G=Z_4(G)$ normalize $E$. Hence $E\cap Z(G)\ne 1$, which is a contradiction.\qed

\begin{thm}{\rm\cite[Theorem 1.5]{GHM}}\label{thm1}
Let $G$ be a $p$-group and $V$ a faithful $\mathbb{F}_p[G]$-module such that there is
no quadratic element in $\Omega_1(Z(G))$.
\begin{itemize*}
  \rr1 If $A$ is an abelian normal subgroup of $G$, then $A$ does not contain any offender.
  \rr2 Suppose that $E$ is an offender. Then $[G', E]\ne 1$.
\end{itemize*}
\end{thm}

\begin{thm}\label{EK}
Let $G$ be a $p$-group and $V$ a faithful $\mathbb{F}_p[G]$-module. If there is a quadratic offender $E$ such that $E^G=EK$ where $K$ is an abelian normal subgroup of $G$, then there is a quadratic element in $\Omega_1(Z(G))$.
\end{thm}
\demo Let $K$ be a maximal member such that $E^G=EK$ where $K$ is an abelian normal subgroup of $G$.

If $K=E^G$, then $K$ contains a offender $E$. By Theorem \ref{thm1} (1), there is a quadratic element in $\Omega_1(Z(G))$.

If $K\ne E^G$, then there is an $a\in E$ such that $N=\lg a\rg K\unlhd G$. Hence $N'\cap \Omega_1(Z(G))\ne 1$. Since $K$ is abelian, $N'=\{[a,x]\mid x\in K\}$. Hence there exists an $x\in K$ such that $1\ne [a,x]\in \Omega_1(Z(G))$. By Lemma \ref{quadratic element}, $[a,x]$ is quadratic.
\qed

\section{$p$-groups of maximal class}

Firstly we recall some terminology. Let $G$ be a
$p$-group. The (characteristic) subgroup of $G$ generated by $\{g^p \mid g\in G\}$ is denoted
by $\mho_1(G)$. Now suppose that $G$ has order $p^n$ and is of maximal class. Set
$G_1 = C_G(G_2/G_4)$. One says that $G$ is {\it exceptional} if $n\ge 5$ and there is some
$3\le i \le n-2$ such that $C_G(G_i/G_{i+2})\ne G_1$.

\begin{defi}{\rm\cite[Definition 3.2.1]{LM}} The degree of commutativity $l$ of a $p$-group $G$ of maximal class is defined to be the maximum integer such that $[G_i,G_j]\le G_{i+j+l}$ for all $i,j\ge 1$ if $G_1$ is not abelian, and $l=n-3$ if $G_1$ is abelian.
\end{defi}

\begin{thm} {\rm\cite[Corollary 3.2.7]{LM}}\label{1} The degree of commutativity of $G$ is positive if and only if $G$ is not exceptional.
\end{thm}

\begin{thm} {\rm\cite[Corollary 3.3.4 (i)]{LM}}\label{2} Let $G$ be a $p$-group of maximal class of order $p^n$ with $n>p+1$. Then $\Omega_1(G_1)=G_{n-p+1}$.
\end{thm}

\begin{thm} {\rm\cite[Theorem 3.3.5]{LM}}\label{3} Let $G$ be a $p$-group of maximal class of order $p^n$ with $n>p+1$. Then $G$ has positive degree of commutativity.
\end{thm}

\begin{thm} {\rm\cite[Theorem 3.2.11]{LM}}\label{4} Let $G$ be a $p$-group of maximal class of order $p^n$ where $n$ is odd and $5\le n\le 2p+1$. Then $G$ has positive degree of commutativity.
\end{thm}

\begin{thm}\label{maximal class}
Let $p$ be an odd prime and $G$ a $p$-group of maximal class of order $p^n$. If $n\ge \max\{2p-6,p+2\}$ then $G$ satisfies
Quadratic Conjecture.
\end{thm}
\demo Otherwise, there is an F-module $V$ such that there is
no quadratic element in $\Omega_1(Z(G))$. By Theorem \ref{weakly}, there is a weakly closed quadratic offender $E$ on $V$.
If $E\trianglelefteq G$, then $G_{n-1}=Z(G)\le E$ and hence there is a quadratic element in $\Omega_1(Z(G))$. In the following, we may assume that $E<E^G$.
By Theorem \ref{thm1} (1), we may assume that $E^G$ is not abelian.

By Lemma \ref{Z_i(G)} (3), $[Z_3(G),E]=1$. Hence $E\le C_G(G_{n-3}/G_{n-1})$. Since $n\ge p+2$, by Theorem \ref{1} $\&$ \ref{3}, $G$ is not exceptional. Hence $E\le G_1$. Furthermore, $E\le \Omega_1(G_1)$. By Theorem \ref{2}, $\Omega_1(G_1)=G_{n-p+1}$. It follows that $E\le G_{n-p+1}$.

Since $E^G$ is not abelian, $1\ne [E,E^G]$. Since $G$ has positive degree of commutativity, $[G_{n-p+1},G_{n-p+1}]=[G_{n-p+1},G_{n-p+2}]\le G_{2n-2p+4}$. Since $n\ge 2p-6$, $2n-2p+4\ge n-2$. Hence $1\ne [E,E^G]\le [G_{n-p+1},G_{n-p+1}]\le G_{n-2}=Z_2(G)$. Thus there are $a\in E$ and $b\in E^G$ such that $1\ne c:=[a,b]\in Z_2(G)$. By Lemma \ref{Z_i(G)} (2), $c\in C_G(a,b)$. By Lemma \ref{quadratic element}, $c$ is a quadratic element. This contradicts Lemma \ref{Z_i(G)} (1).
\qed

\medskip

By Theorem \ref{maximal class}, if $G$ is a $5$-group of maximal class of order $5^n$ with $n\ge 7$, then $G$ satisfies
Quadratic Conjecture; and if $G$ is a $7$-group of maximal class of order $7^n$ with $n\ge 9$, then $G$ satisfies
Quadratic Conjecture. By Corollary \ref{class 4}, if $G$ is a $p$-group of maximal class of order $p^n$ with $n\le 5$, then $G$ satisfies
Quadratic Conjecture. In the following, we prove that if $G$ is a $p$-group of maximal class of order $p^n$ with $6\le n\le 8$, then $G$ satisfies
Quadratic Conjecture. It will deduce that Quadratic Conjecture holds if $G$ is a $p$-group of maximal class where $p\le 7$.

\begin{thm}
Let $G$ be a $p$-group of maximal class with $|G|=p^6$. Then $G$ satisfies
Quadratic Conjecture.
\end{thm}
\demo Otherwise, there is an F-module $V$ such that there is
no quadratic element in $\Omega_1(Z(G))$. By Theorem \ref{weakly}, there is a weakly closed quadratic offender $E$ on $V$. By Lemma \ref{Z_i(G)} (4), $G'=Z_4(G)$ normalize $E$ and $[G',E]\le G_3\cap E$.
 By Theorem \ref{thm1} (2), $[G',E]\ne 1$. Hence $E\cap G_3\ne 1$. Since there is
no quadratic element in $Z(G)=G_5$, $[E\cap G_3,G']\le E\cap G_5=1$. It follows that $E\cap G_3\le Z(G')$. By \cite[Theorem 6.2]{GHM}, $G'$ is not abelian. Hence $Z(G')\le G_4=Z_2(G)$. Thus $E\cap G_3\le Z_2(G)$. This contradicts Lemma \ref{Z_i(G)} (2).
\qed

\begin{thm}
Let $p$ be an odd prime and $G$ a $p$-group of maximal class with $|G|=p^7$. Then $G$ satisfies
Quadratic Conjecture.
\end{thm}
\demo Otherwise, there is an F-module $V$ such that there is
no quadratic element in $\Omega_1(Z(G))$. By Theorem \ref{weakly}, there is a weakly closed quadratic offender $E$ on $V$. By Lemma \ref{Z_i(G)} (4), $G_3=Z_4(G)$ normalize $E$.

By Lemma \ref{Z_i(G)} (3), $[Z_3(G),E]=1$. Hence $E\le C_G(G_{4}/G_{6})$. By Theorem \ref{1} $\&$ \ref{4}, $G$ is not exceptional. Hence $E\le G_1$. It follows that $[G_3,E]\le E\cap [G_3,G_1]\le E\cap G_5=1$.

 Since $[G_3,G_3]\le G_7=1$, $G_3$ is abelian. By Theorem \ref{thm1} (1), $E\not\le G_3$.
Hence $G'\le E^{G}\le G_{1}$.

Since $[E,G_{3}]=1$, $[E^{G},G_{3}]=1$. It follows that $[G',G']=[G',G_{3}]\le [E^G,G_3]=1$. Hence $G'$ is abelian. This contradicts \cite[Theorem 6.2]{GHM}.\qed

\begin{lem}\label{tec}
 Suppose that $A$ is an abelian normal subgroup of a non-abelian $p$-group $G$, $G/A=\lg xA\rg$ is cyclic, $B\le A$ and $B\unlhd G$. Then $[B,G]\cong B/B\cap Z(G)$.
\end{lem}
\demo Since $G/A$ is cyclic, $G'\le A$. Let $\varphi: B\mapsto [B,G]$ be a mapping defined as follows: $\varphi(b)=[b,x]$ for $b\in B$. If $b_1,b_2\in B$, then
$$\varphi(b_1b_2)=[b_1b_2,x]=[b_1,x]^{b_2}[b_2,x]=[b_1,x][b_2,x]=\varphi(b_1)\varphi(b_2),$$
since $[b_1,x], {b_2}\in A$ and $A$ is abelian. Hence $\varphi$ is a homomorphism. If $b\in B$ and $g\in G$, Then
  $$[b,x]^g=[b^g,x^g]=[b^g,x[x,g]]=[b^g,[x,g]][b^g,x]^{[x,g]}=[b^g,x],$$
since $b^g,[x,g],[b^g,x]\in A$ and $A$ is abelian. It follows that $\Im\varphi$ is normal in $G$. Since $B$ centralizes $G/\Im\varphi$, $[B,G]\le \Im\varphi$. Hence $\Im\varphi=[B,G]$. Next, $\Ker\varphi=\{ b\in B\mid [b,x]=1\}=C_B(x)=B\cap Z(G)$. Hence $\Im\varphi=[B,G]\cong B/B\cap Z(G)$.
\qed

\begin{thm}
Let $p$ be an odd prime and $G$ a $p$-group of maximal class with $|G|=p^8$. Then $G$ satisfies
Quadratic Conjecture.
\end{thm}
\demo Otherwise, there is an F-module $V$ such that there is
no quadratic element in $\Omega_1(Z(G))$. By Theorem \ref{weakly}, there is a weakly closed quadratic offender $E$ on $V$. By Lemma \ref{Z_i(G)} (4),
$G_4=Z_4(G)$ normalize $E$.

For convenience, let $G^*=C_G(G_{5}/G_{7})$, a maximal subgroup of $G$.  By Lemma \ref{Z_i(G)} (3), $[G_5,E]=[Z_3(G),E]=1$. Hence
\begin{equation}
E\le C_G(G_{5}/G_{7})=G^*.
\end{equation}
By \cite[Theorem 6.2]{GHM}, $G'$ is not abelian. Hence
\begin{equation}
Z(G')\le G_4.
\end{equation}

We claim that $E\not\le G_3$. Otherwise, $E\le G_3$. In this case, $[E,G_3]\le [G_3,G_3]=[G_3,G_4]\le G_7$. Since $[G_3,E,E]=1$, $G_3$ normalize $E$. Hence $[E,G_3]\le E\cap G_7=1$. It follows that $E\le Z(G_3)$. This contradicts Theorem \ref{thm1} (1).

We claim that $E\not\le G'$. Otherwise, $E^G=G'$. In this case, $[E,G_4]\le [G',G_4]\le G_6$. Since $G_4$ normalize $E$, $[E,G_4]\le E\cap G_6=1$. It follows that $[E^G,G_4]=1$. Hence $G_4=Z(G')$. Note that $G_3$ is abelian and $E^G=EG_3$. By Theorem \ref{EK}, there is
a quadratic element in $\Omega_1(Z(G))$. This contradicts the hypothesis.

By above argument,
\begin{equation}\label{4.3}
E^G=G^*=EG'.
\end{equation}
Since $G'$ is not abelian, $[E^G,G_3]\ge [G',G_3]=G''\ne1$. Hence
\begin{equation}\label{4.4}
[E,G_3]\neq1.
\end{equation}
{\bf Case 1.} $[E\cap G', G']\ne 1$.

In this case, we claim that
\begin{equation}\label{4.5}
E\cap G'\not\le G_4.
\end{equation}
Otherwise, $1\ne [E\cap G', G']\le [G_4,G']\le G_6$. Hence there are $a\in E\cap G'$ and $b\in G'$ such that $1\ne c:=[a,b]\in G_6=Z_2(G)$. Since $c\in C_G(a,b)$, by Lemma \ref{quadratic element}, $c$ is a quadratic element. This contradicts Lemma \ref{Z_i(G)} (1).

Since $G_4$ normalize $E$, $[E\cap G', G_4]\le G_6\cap E$.  By Lemma \ref{Z_i(G)} (2), $G_6\cap E=Z_2(G)\cap E=1$. Hence $[E\cap G',G_4]=1$. Furthermore,
\begin{equation}\label{4.6}
[(E\cap G')^G, G_4]=1.
\end{equation}

We claim that $E\cap G'\le G_3$. Otherwise, $(E\cap G')^G=G'$. By (\ref{4.6}), $G_4\le Z(G')$. It follows that $G_3$ is abelian. By Lemma \ref{tec}, $G''=[G',G_3]\cong G_3/G_3\cap Z(G')=G_3/G_4$. Hence $G''=G_7$. Since $1\ne [E\cap G',G']=G_7=Z(G)$, there are $a\in E\cap G'$ and $b\in G'$ such that $1\ne c:=[a,b]\in Z(G)$. By Lemma \ref{quadratic element}, $c$ is a quadratic element. This contradicts the hypothesis. By (\ref{4.5}),
\begin{equation}\label{4.7}
(E\cap G')^G=G_3.
\end{equation}
By (\ref{4.6}) and (\ref{4.7}), $[G_3,G_3]=[G_3,G_4]=[(E\cap G')^G,G_4]=1$. Hence $G_3$ is abelian. By Lemma \ref{tec}, $[G',G_4]\cong G_4/G_4\cap Z(G')=G_4/Z(G')$. Since $[G_5,E]=1$, $[G_5,E^G]=1$. Hence $G_5\le Z(G')$ and $|[G',G_4]|\le p$. Thus
\begin{equation}
[G',G_4]\le G_7.
\end{equation}
Let $K=EG_3$. Then $K/G_3\cong E/E\cap G_3=E/E\cap G'\cong EG'/G'=G^*/G'$ is cyclic. By Lemma \ref{tec}, $[K,G_4]\cong G_4/G_4\cap Z(K)$. Since $G_5\le Z(K)$, $|[K,G_4]|\le p$. It follows that $[G^*,G_4]=[KG',G_4]=[K,G_4][G',G_4]$ is of order at most $p^2$. Hence
\begin{equation}
 [E,G_4]\le [G^*,G_4]\le G_6.
\end{equation}
Furthermore, $[E,G_4]\le G_6\cap E=1$. By (\ref{4.3}), $G_4\le Z(G^*)\le Z(G')$. By Lemma \ref{tec}, $G''=[G',G_3]\cong G_3/G_3\cap Z(G')=G_3/G_4$. Hence $G''=G_7$. Furthermore, $1\ne [E\cap G',G']=G''=G_7=Z(G)$. Hence there are $a\in E\cap G'$ and $b\in G'$ such that $1\ne c:=[a,b]\in Z(G)$. By Lemma \ref{quadratic element}, $c$ is a quadratic element. This contradicts the hypothesis.

\noindent {\bf Case 2.} $[E\cap G', G']=1$.

In this case, $E\cap G'\le Z(G^*)$ since $[E\cap G',E]=1$. Hence $E\cap G'\le G_4$ and
\begin{equation}
E\cap G'=E\cap G_4.
\end{equation}
By Lemma \ref{Z_i(G)} (2), $E\cap G_6=1$. It follows that $|E\cap G_4|\le p^2$. On the other hand, by Lemma \ref{p3}, $|E|\ge p^3$. It follows that $|E\cap G_4|=|E\cap G'|\ge p^2$. Thus $|E|=p^3$, $|E\cap G_4|=p^2$, $|E\cap G_5|=p$ and $E\cap G'\not\le G_5$. Hence
\begin{equation}
Z(G^*)=G_4=Z(G'),
\end{equation}
and $G_3$ is abelian. Note that $EG_3/G_3\cong E/E\cap G_3=E/E\cap G'\cong EG'/G'$ is cyclic. By Lemma \ref{tec}, $[EG_3,G_3]\cong G_3/G_3\cap Z(EG_3)=G_3/G_4$, and $[G',G_3]\cong G_3/G_3\cap Z(G')=G_3/G_4$. Hence $[E,G_3]$ and $[G',G_3]$ are of order $p$. It follows that $[G^*,G_3]=[EG',G_3]=[E,G_3][G',G_3]\le G_6$. In this case, $[G_3,E,E]=1$. Hence $G_3$ normalize $E$. Thus $[E,G_3]\le E\cap G_6=1$. This contradicts (\ref{4.4}).\qed

%{\bf Acknowledgments}  I cordially thank the referee for detail
%reading and helpful comments, which helped me to improve the whole paper considerably.


\begin{thebibliography}{99}

%\bibitem{CD}
% A. Chermak, A. Delgado,  A measuring argument for finite groups, {\it Proc. Amer. Math. Soc.} 107(1989), no.4, 907--914.
%
%\bibitem{C}
%Chermak, A.: Quadratic action and the $\mathscr{P}(G, V)$-theorem in arbitrary characteristic. J. Group Theory
%2(1), 1-13 (1999)

\bibitem{GHL}
D. J. Green, L. H$\acute{e}$thelyi and M. Lilienthal, On Oliver's $p$-group conjecture, {\it Alg.
Number Theory} {\bf 2} (2008), 969-977.

\bibitem{GHM}
D. J. Green, L. H$\acute{e}$thelyi and N. Mazza, On Oliver's $p$-group conjecture: II, {\it Math.
Ann.} {\bf 347} (2010), 111-122.

\bibitem{GHM2}
D. J. Green, L. H$\acute{e}$thelyi and N. Mazza, On a strong form of Oliver's $p$-group conjecture,
J. Algebra, {\bf 342} (2011), 1-15.

\bibitem{GLS}
Gorenstein, D., Lyons, R., Solomon, R.: The classification of the finite simple groups. Number 2,
Mathematical Surveys and Monographs, vol. 40. American Mathematical Society, Providence, RI
(1996)

%\bibitem{Huppert}
%B. Huppert, Endliche Gruppen, I. Die Grundlehren der Mathematischen Wissenschaften, Band 134.
%Springer, Berlin (1967)
%
%\bibitem{HB}
%B. Huppert and N. Blackburn, Finite groups, III (Springer, 1982).

\bibitem{LM}
Leedham-Green,C.R., McKay, S.: The structure of groups of prime power order, London Mathematical
Society Monographs. New Series, vol. 27. Oxford University Press, Oxford (2002)

\bibitem{MS}
U. Meierfrankenfeld and B. Stellmacher, The other $\mathscr{P}(G, V)$-theorem, {\it Rend. Sem. Mat. Univ. Padova}
{\bf 115}, (2006), 41-50.

\bibitem{O}
B. Oliver, Equivalences of classifying spaces completed at odd primes, {\it Math. Proc. Cambridge Philos. Soc.} {\bf 137}(2), (2004),
321-347.

\end{thebibliography}
\end{document}